\documentclass[10pt]{amsart}
\usepackage{xcolor}
\usepackage{hyperref}
\hypersetup{
	colorlinks=true,
	citecolor=blue,
	linkcolor=blue,
	filecolor=magenta,      
	urlcolor=cyan,
}
\usepackage[capitalise,noabbrev]{cleveref}
\usepackage[all, knot]{xy}\SelectTips{eu}{}

\usepackage[all]{xy}
\usepackage{amsmath,amscd}
\usepackage{amssymb}
\usepackage{mathtools}
\usepackage{stmaryrd}
\usepackage{url}
\usepackage{tikz-cd}
\usepackage{enumitem}

\usepackage{listings}
\usepackage[breakable, listings]{tcolorbox}

\lstdefinelanguage{lean}{
  keywords={def,theorem,lemma,structure,class,where,fun,let,in,do,
             match,with,if,then,else,return,have,show,from,by,
             exists,forall,noncomputable,instance,abbrev,variable,
             open,namespace,end,import,universe,Sort,Prop,Type},
  sensitive=true,
  morecomment=[l]{--},
  morecomment=[s]{/-}{-/},
  morestring=[b]",
}

\newtcblisting{lean}{
  colback=gray!5,
  colframe=gray!20,
  sharp corners,
  boxrule=0.5pt,
  breakable=false,
  listing only,
  listing options={language=lean, basicstyle=\small\ttfamily, breaklines=true},
}

\usepackage{setspace}
\usepackage{microtype}

\usepackage{color}

\usepackage{tikz}

\usepackage{enumitem,kantlipsum}

\usepackage{unicode-math}
\definecolor{leanblue}{rgb}{0.18, 0.33, 0.59}
\definecolor{leangreen}{rgb}{0.13, 0.55, 0.13}
\definecolor{leangray}{rgb}{0.45, 0.45, 0.45}
\definecolor{leanpurple}{rgb}{0.5, 0.0, 0.5}

\newtheorem{thm}{Theorem}[section]
\newtheorem*{thm*}{Theorem}
\newtheorem*{cor*}{Corollary}
\newtheorem*{dfn*}{Definition}

\newtheorem{cthm}{Theorem}

\newtheorem{lem}[equation]{Lemma}

\theoremstyle{definition}

\newcommand{\inj}{\hookrightarrow}

\title{Formalization in Lean of faithfully flat descent of projectivity}
\author{Liran Shaul}
\address{Department of Algebra, Faculty of Mathematics and Physics, Charles University in Prague, Sokolovsk\'a 83, 186 75 Praha, Czech Republic}

\email{shaul@karlin.mff.cuni.cz}

\thanks{{\em Mathematics Subject Classification} 2020:
13C10, 16D40, 68V20}

\begin{document}
\begin{abstract}
We formalize in Lean the following foundational result in commutative algebra: Let $R \to S$ be a faithfully flat map of (not necessarily noetherian) commutative rings, and let $P$ be an arbitrary $R$-module. 
Then $P$ is projective over $R$ if and only if $S\otimes_R P$ is projective over $S$. 
This formalizes and verifies Perry's fix of a subtle gap in the classical work of Raynaud and Gruson, a result which is a key ingredient in the study of finitistic dimension of commutative noetherian rings.
\end{abstract}

\numberwithin{equation}{section}
\maketitle

\setcounter{section}{-1}
\section{Introduction}

The foundational paper \cite{RG} has been one of the most influential papers of 20th-century  commutative algebra. 
One of its major achievements, 
solving a problem posed by Bass in \cite{Bass}, 
was the result stating that the (big projective) finitistic dimension of any commutative noetherian ring is equal to its Krull dimension.
This can be stated, in less technical language,
as the fact that over a commutative noetherian ring $A$,
any module that can be expressed using finitely many (possibly infinite sized) matrices,
can be expressed using at most $\dim(A)$ number of matrices.

To prove this result, as is customary in many problems in modern commutative algebra,
the authors proceed in three steps. 
First, the result is proved for complete noetherian local rings,
then, it is deduced for arbitrary noetherian local rings using faithful flatness,
and finally, a local-to-global argument establishes the result for all commutative noetherian rings. All three steps are highly non-trivial for this problem. 

The passage from complete local rings to arbitrary noetherian local rings is contained in:

\begin{cthm}\label{thm:main}
Let $R \to S$ be a faithfully flat map of commutative rings, and let $P$ be an $R$-module.
Then $P$ is projective over $R$ if and only if the $S$-module $S\otimes_R P$ is projective over $S$.
\end{cthm}

Note that it is not assumed that the rings are noetherian,
and that $P$ could be infinitely generated.

As later observed by Gruson \cite{Gruson},
the proof of \cref{thm:main} in \cite{RG} contained a gap.

A complete proof of \cref{thm:main} was given in \cite[Theorem 9.6]{perry},
and was also incorporated into the stacks project in \cite[Tag 05A9]{SP},
but to our knowledge, it was never peer reviewed.

The aim of this paper is to describe a recent successful formalization of this result,
as well as all its prerequisites in Lean.
All Lean code is in \cite{code}.

The formalization of \cref{thm:main} required a systematic development of pushouts in the category of modules, universal injectivity via domination criteria, a fully formal treatment of the Mittag–Leffler condition for inverse systems, and a transfinite devissage framework using ordinals. All these components were not previously available in Lean 4 and were developed from scratch, resulting in over 10,000 lines of code.

The notion of a Mittag-Leffler module was first defined and studied in \cite{RG}.
This notion, its basic properties, 
and its various applications is formalized in this project in great length.
Perhaps the most interesting result we formalize, other than \cref{thm:main},
is the following result which is precisely \cite[Theorem 7.4]{perry}:

\begin{cthm}\label{thm:projchar}
Let $R$ be a ring, and let $P$ be a left $R$-module.
Then $P$ is projective over $R$ if and only if the following holds:
\begin{enumerate}
\item $P$ is flat over $R$.
\item $P$ is Mittag-Leffler over $R$.
\item $P$ is a direct sum of countably generated $R$-modules.
\end{enumerate}
\end{cthm}

Here, the ring $R$ is possibly noncommutative.
The fact that projective modules satisfy (3) is a classical result of Kaplansky which we also formalize. 
The proof also depends on Lazard's theorem,
which states that any flat module is a direct limit of finitely generated free modules,
and that result, too, is formalized.

In the remainder of this document we describe the most important results we formalize,
including their exact statement in Lean,
as well as some of the formalization challenges.
All results described here were successfully compiled using Lean 4 with Mathlib, and we intend to submit them for inclusion in Mathlib.

\textbf{Acknowledgments.}

The author thanks Mirek Olsak as well as the Zulip community for helpful discussions.
The author acknowledges the assistance of the large language models Claude Opus 4.5 and Claude Opus 4.6 in generating some of the Lean code used in this project.
This work has been supported by the grant GA~\v{C}R 26-22734S from the Czech Science Foundation.

\section{Kaplansky Devissage}

The results of this section are in the file lib/Kap.lean in \cite{code}.
A key element in the proof of both \cref{thm:main} and \cref{thm:projchar} is the reduction to the case of countably generated modules.
In order to be able to do that, we require the following result due to Kaplansky \cite[Theorem 1]{Kaplansky}, which we formalized:

\begin{thm}\label{thm:kap}
Let $R$ be a ring, 
let $M$ be an $R$-module,
and suppose that $M$ is a direct sum of countably generated $R$-modules.
If $P$ is an $R$-module which is a direct summand of $M$,
then $P$ is also a direct sum of countably generated $R$-modules.
\end{thm}

In order to formalize this result in Lean,
we first define a module to be countably generated if it has a countable subset which generates it:
\begin{lean}
def Submodule.IsCountablyGenerated
  {R : Type u} {M : Type v}
  [Ring R] [AddCommGroup M] [Module R M]
  (N : Submodule R M) : Prop :=
  ∃ s : Set N, s.Countable ∧ span R s = ⊤
\end{lean}

To avoid unnecessary complicated use of quotients modules,
which are not strictly needed here as we only work with direct summands, 
we introduce the following auxiliary definition of direct sum complement

\begin{lean}
def IsRelativeComplement
  {R : Type u} {M : Type v}
  [Ring R] [AddCommGroup M] [Module R M]
  (A B C : Submodule R M) : Prop :=
  A ≤ B ∧ C ≤ B ∧ Disjoint A C ∧ A ⊔ C = B
\end{lean}

Using these definitions, 
we then follow Kaplansky and introduce the notion of a Kaplansky devissage,
given by:

\begin{lean}
structure KaplanskyDevissage
    (R : Type u) [Ring R]
    (M : Type v) [AddCommGroup M] [Module R M]
    (S : Ordinal.{w}) where
  seq : Set.Iio S → Submodule R M
  monotone : ∀ α β : Set.Iio S, α ≤ β → seq α ≤ seq β
  zero_eq_bot (h : 0 < S) : seq ⟨0, h⟩ = ⊥
  union_eq_top : (⨆ α : Set.Iio S, seq α) = ⊤
  limit_continuity : ∀ α : Set.Iio S,
    Order.IsSuccLimit α.val →
    seq α = ⨆ β : {x : Set.Iio S // x.val < α.val}, seq β.val
  succ_step : ∀ (α : Ordinal.{w}) (h : Order.succ α < S),
    ∃ C : Submodule R M,
      IsRelativeComplement
        (seq ⟨α, (Order.lt_succ α).trans h⟩)
        (seq ⟨Order.succ α, h⟩)
        C ∧
      C.IsCountablyGenerated
\end{lean}

Thus, such a structure associated to an $R$-module $M$ is given by some ordinal $S$,
and an $S$-sized filtration of $M$ with the property that each element of it is a direct summand of $M$, 
and each successive quotient in it is countably generated.
The use of 
\begin{lean}
seq : Set.Iio S → Submodule R M
\end{lean}
instead of the possibly more common
\begin{lean}
seq : ∀ α : Ordinal, α < S → Submodule R M
\end{lean}
is to ensure that we stay at the same universe level,
for otherwise Lean, ignoring the $\alpha < S$ condition, 
jumps from the universe $w$ to $w+1$ unnecessarily.

We then prove, following \cite[tag 058V]{SP} and \cite[tag 058W]{SP} that for an $R$-module $M$ over some ring $R$ the existence of a Kaplansky devissage for $M$ with respect to some ordinal $S$ is equivalent to $M$ being a direct sum of countably generated modules:

\begin{lean}
theorem kaplansky_devissage_iff_direct_sum
    {R : Type u} [Ring R] {M : Type v} [AddCommGroup M] [Module R M] :
    (∃ (ι : Type w) (fam : ι → Submodule R M),
      DirectSum.IsInternal fam ∧ ∀ i, (fam i).IsCountablyGenerated)
    ↔
    ∃ (S : Ordinal.{w}), Nonempty (KaplanskyDevissage R M S) 
\end{lean}

The proof formalizes the equivalence by constructing explicit filtrations in both directions, with the harder backward direction requiring a two-level induction, first reducing independence to a finite set of indices via 
\begin{lean}
Submodule.mem_span_finite_of_mem_span    
\end{lean}
then, the contributors are eliminated one by one by peeling off the maximal-ordinal index using the filtration's disjointness conditions.

Using this characterization, we then proceed to prove \cref{thm:kap}.
We are given an $R$-module $M = \oplus_{\alpha \in I} C_{\alpha}$, 
which is a direct sum of countably generated modules $C_{\alpha}$,
and a direct summand $N$ of $M$.
Suppose $N \oplus K = M$. 
The goal is to show $N$ is a direct sum of countably generated $R$-modules.
Using
\begin{lean}
theorem kaplansky_devissage_iff_direct_sum
\end{lean}
above, it is sufficient to construct a Kaplansky devissage for $N$.
Following sources like \cite[Theorem 1]{Kaplansky} and \cite[tag  058X]{SP},
we do this by producing a Kaplansky devissage for $M$ with the extra property that for each $\alpha$ it holds that 
\begin{equation}\label{eqn:decomp}
M_{\alpha} = (M_{\alpha} \cap N) \oplus (M_{\alpha} \cap K)
\end{equation}
To do that, choose some well ordering of the index set $I$.
Assuming an $M_{\alpha}$ was defined, we produce $M_{\alpha+1}$ by 
first taking $j \in I$ to be the smallest $j$ such that $C_j \nsubseteq M_{\alpha}$.
We then define $M_{\alpha+1}$ to be the module generated by $M_{\alpha}$,
the generators of $C_j$,
and throw in generators of other indices from $I$,
so that a decomposition of the form \cref{eqn:decomp} remains true for $M_{\alpha+1}$.
This requires only countably many indices,
which ensures that $M_{\alpha+1}/M_{\alpha}$ is countably generated.
Following this approach we obtain \cref{thm:kap} in the form:
\begin{lean}
theorem dirsummand_of_dirSumCountable {R : Type u} [Ring R]
    {M N : Type v} [AddCommGroup M] [AddCommGroup N] [Module R M] 
    [Module R N] (exDecom : ∃ (ι : Type w) (fam : ι → Submodule R M),
      DirectSum.IsInternal fam ∧ ∀ i, (fam i).IsCountablyGenerated)
    (i : N →ₗ[R] M) (π : M →ₗ[R] N) (hπi : π.comp i = LinearMap.id) :
    ∃ (κ : Type w) (fam' : κ → Submodule R N),
      DirectSum.IsInternal fam' ∧ ∀ j, (fam' j).IsCountablyGenerated    
\end{lean}

Since a projective module is a direct summand of a free module,
and a free module is a direct sum of cyclic modules,
this result immediately implies that any projective module is a direct sum of countably generated projective modules. We thus obtain:
\begin{lean}
theorem proj_is_dirSumCountable {R : Type u} {P : Type u} [Ring R] 
[AddCommGroup P] [Module R P] (pproj : Module.Projective R P) 
: IsDirectSumOfCountablyGenerated R P    
\end{lean}

\section{Lazard theorem}

The results of this section are in the file lib/Lazard.lean in \cite{code}.
Following the results of the previous section, 
the proof of \cref{thm:projchar} then reduces immediately to the problem of showing a countably generated flat module which is Mittag-Leffler is projective.
Proving this requires the classical Lazard theorem \cite{Lazard}:
\begin{thm}\label{thm:Lazard}
Let $R$ be a ring, and let $F$ be an $R$-module.
Then $F$ is flat over $R$ if and only if $F$ is isomorphic to a direct limit of finitely presented free $R$-modules.
\end{thm}

Closely related to this, 
and also crucial for working with Mittag-Leffler,
is an even more basic fact, which was also missing from Mathlib,
namely, that any module over any ring is a direct limit of finitely presented modules.

We start by proving this more basic result,
formulated as:
\begin{lean}
theorem Module.isDirectLimit_of_finitelyPresented :
    ∃ (ι : Type u) (_ : Preorder ι) (_ : IsDirected ι (· ≤ ·)) 
    (_ : Nonempty ι) (_ : DecidableEq ι) (G : ι → Type u) 
    (_ : ∀ i, AddCommGroup (G i)) (_ : ∀ i, Module R (G i))
      (_ : ∀ i, Module.FinitePresentation R (G i))
      (f : ⦃i j : ι⦄ → i ≤ j → G i →ₗ[R] G j) 
      (_ : DirectedSystem G (fun _ _ h => f h)),
      Nonempty (DirectLimit G (fun _ _ h => f h) ≃ₗ[R] M)    
\end{lean}

We then proceed to prove Lazard's theorem.
A key technical lemma needed for this is:
\begin{lem}
Let $R$ be a ring and $M$ be a flat left $R$-module. Let $J \subset M \times \mathbb{Z}$ be a finite set, and let $N \subseteq R^{(J)}$ be a finitely generated $R$-submodule. Let $\psi_J: R^{(J)} \to M$ be the $R$-linear map defined by $e_{(m, z)} \mapsto m$. 

If $N \subseteq \ker(\psi_J)$, then there exists a finite set $J'$ containing $J$ and a finitely generated submodule $N' \subseteq R^{(J')}$ such that:
\begin{enumerate}
    \item $N' \subseteq \ker(\psi_{J'})$.
    \item The quotient $R^{(J')} / N'$ is a free $R$-module.
    \item Under the natural inclusion $\iota: R^{(J)} \hookrightarrow R^{(J')}$, we have $\iota(N) \subseteq N'$.
\end{enumerate}
\end{lem}

which we prove as:
\begin{lean}
lemma Module.Flat.enlarge_to_free [Module.Flat R M]
    {J : Finset (M × ℤ)} {N : Submodule R (J →₀ R)} (hN : N.FG)
    (hN_ker : ∀ x ∈ N, Finsupp.embDomain ⟨Subtype.val, 
    Subtype.val_injective⟩ x ∈ 
    LinearMap.ker (Finsupp.linearCombination R (Prod.fst : M × ℤ → M))) :
    ∃ (J' : Finset (M × ℤ)) (hJJ' : J ⊆ J') (N' : Submodule R (J' →₀ R)) 
    (_ : N'.FG), (∀ x ∈ N', Finsupp.embDomain ⟨Subtype.val, 
    Subtype.val_injective⟩ x ∈
        LinearMap.ker (Finsupp.linearCombination R (Prod.fst : M × ℤ → M))) 
        ∧ Module.Free R ((J' →₀ R) ⧸ N') ∧
      ∀ x ∈ N, Finsupp.lmapDomain R R (fun i => ⟨i.1, hJJ' i.2⟩) x ∈ N'    
\end{lean}

This enables us to construct a directed system with the needed properties,
and obtain \cref{thm:Lazard} in the form of:
\begin{lean}
theorem Module.Flat.Lazard :
    Module.Flat R M ↔
    ∃ (ι : Type u) (_ : Preorder ι) (_ : IsDirected ι (· ≤ ·)) 
    (_ : Nonempty ι) (_ : DecidableEq ι) (G : ι → Type u) 
    (_ : ∀ i, AddCommGroup (G i)) (_ : ∀ i, Module R (G i))
    (_ : ∀ i, Module.FinitePresentation R (G i)) 
    (_ : ∀ i, Module.Free R (G i)) (f : ⦃i j : ι⦄ → i ≤ j → G i →ₗ[R] G j) 
    (_ : DirectedSystem G (fun _ _ h => f h)),
    Nonempty (DirectLimit G (fun _ _ h => f h) ≃ₗ[R] M)    
\end{lean}

\section{Pushout of modules}

The results of this section are in the file lib/Pushout.lean in \cite{code}.
From now on, we specialize to the case where the ring $R$ is commutative. 
If $R$ is a ring, 
and $f:A \to B$ and $g:A \to C$ are two $R$-linear maps between $R$-modules,
their pushout is an $R$-module $D$ and a pair of maps $f':B \to D$ and $g':C \to D$ making the diagram
\[
\xymatrix{
A \ar[r]^f \ar[d]_g & B \ar[d]^{f'}\\
C \ar[r]_{g'} & D
}
\]
commutative, which is universal with respect to this property.
Working in Lean, instead of working with such a universal construction,
it is of course easier to choose a concrete realization of the pushout.
This can be given by taking $D = (B\oplus C)/T$,
where $T$ is the submodule of $B\oplus C$ generated by the set 
\[
\{(f(a),-g(a)) \mid a \in A\}.
\]
The maps $f'$ and $g'$ are then given by composing the inclusions $B \inj B\oplus C$
and $C \inj B\oplus C$ with the natural surjection $B\oplus C \to (B\oplus C)/T$.

We implement this in Lean using the following definitions.
The submodule $T$ is given by
\begin{lean}
def pushoutSubmodule (f : A →ₗ[R] B) (g : A →ₗ[R] C) : Submodule R (B × C) 
:=  Submodule.span R (Set.range (fun a => (f a, -g a)))    
\end{lean}
the pushout is then
\begin{lean}
def Pushout (f : A →ₗ[R] B) (g : A →ₗ[R] C) : Type u :=
    (B × C) ⧸ pushoutSubmodule f g
\end{lean}
And the maps $f'$ and $g'$ are given by
\begin{lean}
def Pushout.inl (f : A →ₗ[R] B) (g : A →ₗ[R] C) : B →ₗ[R] Pushout f g :=
    (Submodule.mkQ _).comp (LinearMap.inl R B C)

def Pushout.inr (f : A →ₗ[R] B) (g : A →ₗ[R] C) : C →ₗ[R] Pushout f g :=
    (Submodule.mkQ _).comp (LinearMap.inr R B C)
\end{lean}

We then prove various basic properties that the pushout satisfy,
including the universal property,
certain functionality,
and the fact it commutes with tensor products.
More precisely, it holds that:
\begin{thm}
Let $R$ be a commutative ring, $S$ an $R$-algebra, and $f: A \to B$, $g: A \to C$ 
be $R$-linear maps. Denote by 
\[
  f_S = \mathrm{id}_S \otimes f : S \otimes_R A \to S \otimes_R B, 
  \qquad 
  g_S = \mathrm{id}_S \otimes g : S \otimes_R A \to S \otimes_R C
\]
the base-changed maps. Then there exists an $S$-linear isomorphism
\[
  \Phi : S \otimes_R \mathrm{Pushout}(f,g) \xrightarrow{\;\sim\;} \mathrm{Pushout}(f_S, g_S).
\]
Moreover, the following diagram commutes:
\[
\xymatrix@C=6em@R=4em{
  S \otimes_R C 
    \ar[r]^{\mathrm{id}_S \otimes \iota_C} 
    \ar[dr]_{\iota_C^S} 
  & S \otimes_R \mathrm{Pushout}(f,g) 
    \ar[d]^{\Phi} \\
  & \mathrm{Pushout}(f_S, g_S)
}
\]
where $\iota_C : C \to \mathrm{Pushout}(f,g)$ and 
$\iota_C^S : S \otimes_R C \to \mathrm{Pushout}(f_S, g_S)$ 
are the canonical right inclusions into the respective pushouts.
\end{thm}
which we implement in Lean using:
\begin{lean}
noncomputable def Pushout.baseChangeEquiv
    {R : Type u} [CommRing R]
    {S : Type u} [CommRing S] [Algebra R S]
    {A B C : Type u} [AddCommGroup A] [AddCommGroup B] [AddCommGroup C]
    [Module R A] [Module R B] [Module R C]
    (f : A →ₗ[R] B) (g : A →ₗ[R] C) :
    S ⊗[R] Pushout f g ≃ₗ[S] 
    Pushout (LinearMap.baseChange S f) (LinearMap.baseChange S g)    
\end{lean}
and
\begin{lean}
lemma Pushout.baseChangeEquiv_inr
    {R : Type u} [CommRing R]
    {S : Type u} [CommRing S] [Algebra R S]
    {A B C : Type u} [AddCommGroup A] [AddCommGroup B] [AddCommGroup C]
    [Module R A] [Module R B] [Module R C]
    (f : A →ₗ[R] B) (g : A →ₗ[R] C) :
    (↑(Pushout.baseChangeEquiv (S := S) f g) ∘ₗ
    LinearMap.baseChange S (Pushout.inr f g)) =
      Pushout.inr (LinearMap.baseChange S f) (LinearMap.baseChange S g)    
\end{lean}

\section{Universally injective maps}

The results of this section are in the file lib/UnivInj.lean in \cite{code}.
Given a commutative ring $R$, 
recall that a map of $R$-modules $f:M \to N$ is called a universally injective map,
if for any $R$-module $Q$,
the map $f\otimes_R \mathrm{id} : M\otimes_R Q \to N\otimes_R Q$ is injective.
Such a map is necessarily injective. In Lean, we define it as:
\begin{lean}
def UniversallyInjective (f : M →ₗ[R] N) : Prop :=
  ∀ (Q : Type u) [AddCommGroup Q] [Module R Q], 
  Function.Injective (f.rTensor Q)
\end{lean}

A key property satisfied by universally injective maps is the following lifting property we need,
which is contained in \cite[tag 058K]{SP}:
\begin{thm}
Let $R$ be a commutative ring, let $f : M \to N$ be a universally injective $R$-linear map,
and let $F, G$ be finitely generated free $R$-modules. 
Suppose there exist $R$-linear maps
$g : F \to M$, $h : G \to N$, and $k : F \to G$ such that the square
\[
\xymatrix@C=4em@R=3em{
  F \ar[r]^{k} \ar[d]_{g} & G \ar[d]^{h} \\
  M \ar[r]_{f}             & N
}
\]
commutes, i.e.\ $h \circ k = f \circ g$. Then there exists an $R$-linear map
$\phi : G \to M$ such that
\[
  \phi \circ k = g.
\]
\end{thm}
We formalize it as:
\begin{lean}
theorem universially_injective_lift_free
    {R : Type u} [CommRing R]
    {M : Type u} {N : Type u}
    [AddCommGroup M] [AddCommGroup N]
    [Module R M] [Module R N]
    {F : Type u} [AddCommGroup F] [Module R F] 
    [Module.Free R F] [Module.Finite R F]
    {G : Type u} [AddCommGroup G] [Module R G] 
    [Module.Free R G] [Module.Finite R G]
    (f : M →ₗ[R] N)
    (univf : UniversallyInjective f)
    (g : F →ₗ[R] M)
    (h : G →ₗ[R] N)
    (k : F →ₗ[R] G)
    (commSquare : h ∘ₗ k = f ∘ₗ g) :
    ∃ φ : G →ₗ[R] M, φ ∘ₗ k = g    
\end{lean}
The key idea of the proof is to choose coordinates using chosen bases of $F$ and $G$,
reducing the lifting problem to finding $z_j \in M$ such that
$x_i = \sum_j a_{ij} z_j$ for all basis vectors $e_i$ of $F$,
where $x_i = g(e_i)$ and $a_{ij}$ are the matrix coefficients of $k$
defined by $k(e_i) = \sum_j a_{ij} b_j$,
and then using universal injectivity to conclude that such $z_j$ exist.

Another key property is the fact this notion descends along faithfully flat maps:
\begin{thm}
Let $R$ be a commutative ring, let $R \to S$ be a faithfully flat map of commutative rings,
and let $f : M \to N$ be an $R$-linear map. If the base-changed map
\[
  f_S = \mathrm{id}_S \otimes f : S \otimes_R M \to S \otimes_R N
\]
is universally injective as an $S$-linear map, then $f$ is universally injective as an
$R$-linear map.
\end{thm}
which we formalize as:
\begin{lean}
theorem UniversallyInjective.of_baseChange_faithfullyFlat
    {R : Type u} [CommRing R]
    {S : Type u} [CommRing S] [Algebra R S]
    [Module.FaithfullyFlat R S]
    {M : Type u} {N : Type u}
    [AddCommGroup M] [AddCommGroup N]
    [Module R M] [Module R N]
    (f : M →ₗ[R] N)
    (hf : UniversallyInjective (f.baseChange S)) :
    UniversallyInjective f
\end{lean}

\section{Mittag-Leffler Inverse Systems}

The results of this section are in the file lib/mlSystem.lean in \cite{code}.
We formalize the classical Mittag-Leffler condition for inverse systems, 
using the Mathlib typeclass \texttt{InverseSystem}, which models a contravariant 
functor from a preordered set $\iota$ to the category of types. 
Concretely, an inverse system consists of a family of types $(F_i)_{i \in \iota}$ 
together with transition maps $f_{ij} : F_j \to F_i$ for $i \leq j$, satisfying 
the identity and composition axioms. 
Following \cite[tag 0594]{SP}, 
we say such a system is \emph{Mittag-Leffler} 
if for each $i \in \iota$, the images $\mathrm{im}(f_{ij})$ stabilize as $j$ increases.
In Lean, we formalize this as:
\begin{lean}
class InverseSystem.IsMittagLeffler (F : ι → Type*) [∀ i, Nonempty (F i)]
    (f : ∀ ⦃i j⦄, i ≤ j → F j → F i) [InverseSystem f] where
  stabilization : ∀ i : ι, ∃ j : ι, ∃ (hij : i ≤ j),
    ∀ k : ι, ∀ (hjk : j ≤ k),
      Set.range (f (hij.trans hjk)) = Set.range (f hij)    
\end{lean}

We prove the following foundational result, which states that a Mittag-Leffler 
inverse system over a countable directed index set has a nonempty inverse limit:
\begin{thm}\label{thm:non-empty}
Let $\iota$ be a countable directed preorder, and let $(F_i, f_{ij})$ be a 
Mittag-Leffler inverse system indexed by $\iota$. Then the inverse limit 
$\varprojlim F_i$ is nonempty.
\end{thm}

which we formalize as:
\begin{lean}
theorem nonempty_inverseLimit_of_countable_mittagLeffler
    {ι : Type*} [Preorder ι] [IsDirected ι (· ≤ ·)] [Nonempty ι]
    [Countable ι]
    {F : ι → Type*} [∀ i, Nonempty (F i)]
    (f : ∀ ⦃i j⦄, i ≤ j → F j → F i) [InverseSystem f]
    [IsMittagLeffler F f] :
    Nonempty (InverseLimit f)
\end{lean}

Here, we use the following concrete definition of the inverse limit:
\begin{lean}
def InverseLimit {ι : Type*} [Preorder ι] {F : ι → Type*}
    (f : ∀ ⦃i j⦄, i ≤ j → F j → F i) : Type _ :=
  { x : ∀ i, F i // ∀ i j (hij : i ≤ j), f hij (x j) = x i }    
\end{lean}
We note that Mathlib's built-in \texttt{InverseSystem.limit} is defined as a limit over elements strictly below a given index, which is unsuitable for our purposes, so we introduce the above definition as a new addition to the library.

The main application of \cref{thm:non-empty} is the following theorem about short exact 
sequences of inverse systems, which is the key ingredient connecting the 
Mittag-Leffler condition to the module theory developed in the next section.
Given a short exact sequence of inverse systems $0 \to (A_i) \to (B_i) \to (C_i) \to 0$ 
of $R$-modules over a countable directed index set, if the system $(A_i)$ is 
Mittag-Leffler, then the induced map on inverse limits $\varprojlim B_i \to \varprojlim C_i$ 
is surjective. A bit more generally, it holds that:

\begin{thm}
Let $\iota$ be a countable directed preorder, and let
\[
A_i \xrightarrow{f_i} B_i \xrightarrow{g_i} C_i \to 0
\]
be a right exact sequence of inverse systems of $R$-modules, meaning each $g_i$ is 
surjective and $\mathrm{im}(f_i) = \ker(g_i)$. If the inverse system $(A_i)$ is 
Mittag-Leffler, then the induced map on inverse limits
\[
\varprojlim B_i \to \varprojlim C_i
\]
is surjective.
\end{thm}

which we formalize as:
\begin{lean}
theorem surjective_limit_of_mittagLeffler_exact
    {ι : Type*} [Preorder ι] [IsDirected ι (· ≤ ·)] [Nonempty ι] 
    [Countable ι] {R : Type*} [Ring R]
    {A B C : ι → Type*}
    [∀ i, AddCommGroup (A i)] [∀ i, AddCommGroup (B i)] 
    [∀ i, AddCommGroup (C i)] [∀ i, Module R (A i)] 
    [∀ i, Module R (B i)] [∀ i, Module R (C i)]
    [∀ i, Nonempty (A i)] [∀ i, Nonempty (B i)] [∀ i, Nonempty (C i)]
    (fA : ∀ ⦃i j⦄, i ≤ j → A j →ₗ[R] A i) 
    [InverseSystem (fun i j h => (fA h))]
    (fB : ∀ ⦃i j⦄, i ≤ j → B j →ₗ[R] B i) 
    [InverseSystem (fun i j h => (fB h))]
    (fC : ∀ ⦃i j⦄, i ≤ j → C j →ₗ[R] C i) 
    [InverseSystem (fun i j h => (fC h))]
    (f : ∀ i, A i →ₗ[R] B i)
    (g : ∀ i, B i →ₗ[R] C i)
    (hg_surj : ∀ i, Function.Surjective (g i))
    (hexact : ∀ i, LinearMap.range (f i) = LinearMap.ker (g i))
    (hf_compat : ∀ i j (h : i ≤ j) (a : A j), fB h (f j a) = f i (fA h a))
    (hg_compat : ∀ i j (h : i ≤ j) (b : B j), fC h (g j b) = g i (fB h b))
    [IsMittagLeffler A (fun i j h => (fA h))] :
    Function.Surjective (InverseLimit.map
      (fun i j h => (fB h))
      (fun i j h => (fC h))
      (fun i => (g i))
      (fun i j h => by
        ext x
        exact hg_compat i j h x))
\end{lean}

\section{Domination}

The results of this section are in the file lib/Domination.lean in \cite{code}.
A key technical notion underlying the theory of Mittag-Leffler modules is that of
domination of linear maps. Given $R$-linear maps $f : M \to N$ and $g : M \to M'$,
we say that $g$ \emph{dominates} $f$ if for every $R$-module $Q$,
\[
\ker(f \otimes \mathrm{id}_Q) \subseteq \ker(g \otimes \mathrm{id}_Q).
\]
In Lean, this is formalized as:
\begin{lean}
def LinearMap.Dominates (f : M →ₗ[R] N) (g : M →ₗ[R] M') : Prop :=
  ∀ (Q : Type u) (_ : AddCommGroup Q) (_ : Module R Q),
    ∀ x : TensorProduct R M Q,
      (TensorProduct.map f LinearMap.id) x = 0 →
        (TensorProduct.map g LinearMap.id) x = 0
\end{lean}
Note that in the Lean code, \texttt{Dominates f g} means that $g$ dominates $f$.
We say $f$ and $g$ \emph{mutually dominate} each other if each dominates the other,
i.e.\ $\ker(f \otimes \mathrm{id}_Q) = \ker(g \otimes \mathrm{id}_Q)$ for every
$R$-module $Q$.

A fundamental characterization of domination, which is used throughout the
Mittag-Leffler theory, is the following result, corresponding to \cite[Lemma 6.10]{perry} and \cite[Tag 059D]{SP}:

\begin{thm}
Let $R$ be a commutative ring, and let $f : M \to N$ and $g : M \to M'$ be
$R$-linear maps. Assume that $N/\mathrm{im}(f)$ is finitely presented. Then
$g$ dominates $f$ if and only if $g$ factors through $f$, i.e.\ there exists an
$R$-linear map $h : N \to M'$ such that $g = h \circ f$.
\end{thm}

which we formalize as:
\begin{lean}
theorem LinearMap.dominates_iff_factors_through (f : M →ₗ[R] N)
    (g : M →ₗ[R] M') [Module.FinitePresentation R (N ⧸ LinearMap.range f)] :
    Dominates f g ↔ ∃ (h : N →ₗ[R] M'), g = h.comp f
\end{lean}

The proof makes essential use of both the pushout construction from Section~3
and the universal injectivity theory from Section~4, via the following intermediate
result connecting all three notions:

\begin{thm}\label{thm:domPushout}
Let $R$ be a commutative ring, and let $f : M \to N$ and $g : M \to M'$ be
$R$-linear maps. Then $g$ dominates $f$ if and only if the canonical map
$\iota_C : M' \to \mathrm{Pushout}(f, g)$ is universally injective.
\end{thm}

which we formalize as:
\begin{lean}
theorem LinearMap.dominates_iff_pushout_inr_universallyInjective
    (f : M →ₗ[R] N) (g : M →ₗ[R] M') :
    Dominates f g ↔ UniversallyInjective (Pushout.inr f g)
\end{lean}

\section{Mittag-Leffler modules}

The results of this section are in the file lib/mlModule.lean in \cite{code}.
We now turn to the notion of a Mittag-Leffler module, which was first introduced
by Raynaud and Gruson in \cite{RG}. We continue to work over a commutative ring
$R$. A key role is played by the domination notion from the previous section.
An $R$-module $M$ is called \emph{Mittag-Leffler} if for every finitely presented
$R$-module $P$ and every $R$-linear map $f : P \to M$, there exists a finitely
presented $R$-module $Q$ and an $R$-linear map $g : P \to Q$ such that $f$ and
$g$ mutually dominate each other. In Lean, this is formalized as:
\begin{lean}
def Module.IsMittagLeffler (R : Type u) [CommRing R]
    (M : Type u) [AddCommGroup M] [Module R M] : Prop :=
  ∀ (P : Type u) (_ : AddCommGroup P) (_ : Module R P)
    (_ : Module.FinitePresentation R P),
    ∀ f : P →ₗ[R] M,
      ∃ (Q : Type u) (_ : AddCommGroup Q) (_ : Module R Q)
        (_ : Module.FinitePresentation R Q),
        ∃ g : P →ₗ[R] Q, LinearMap.MutuallyDominate f g
\end{lean}

Using \texttt{Module.isDirectLimit\_of\_finitelyPresented} proved in Section 2,
any $R$-module $M$ can be written as a directed colimit
$M = \varinjlim_{i \in \iota} M_i$ of finitely presented $R$-modules. 
Following \cite[Proposition 6.11]{perry}, and \cite[Tag 0599]{SP}, the Mittag-Leffler condition admits several 
equivalent reformulations in terms of such a presentation, which we now describe.
Given such a presentation, we say it satisfies condition (2) if for each $i \in \iota$,
there exists $j \geq i$ such that the canonical map $M_i \to M$ is dominated by
the transition map $f_{ij} : M_i \to M_j$. In Lean:
\begin{lean}
def IsMittagLeffler' {ι : Type u} [Preorder ι] [IsDirected ι (· ≤ ·)] 
    [Nonempty ι] [DecidableEq ι]
    (F : ι → Type u) [∀ i, AddCommGroup (F i)] [∀ i, Module R (F i)]
    [∀ i, Module.FinitePresentation R (F i)]
    (f : ∀ i j, i ≤ j → F i →ₗ[R] F j)
    [DirectedSystem F (fun i j h => f i j h)] : Prop :=
  ∀ i : ι, ∃ j : ι, ∃ (hij : i ≤ j),
    LinearMap.Dominates (Module.DirectLimit.of R ι F f i) (f i j hij)
\end{lean}
Condition (3) replaces domination by an explicit factorization: for each $i$, there
exists $j \geq i$ such that for all $k \geq i$, the map $f_{ij}$ factors through
$f_{ik}$. In Lean:
\begin{lean}
def IsMittagLeffler'' {ι : Type u} [Preorder ι] [IsDirected ι (· ≤ ·)]
    [Nonempty ι]
    (F : ι → Type u) [∀ i, AddCommGroup (F i)] [∀ i, Module R (F i)]
    [∀ i, Module.FinitePresentation R (F i)]
    (f : ∀ i j, i ≤ j → F i →ₗ[R] F j) : Prop :=
  ∀ i : ι, ∃ j : ι, ∃ (hij : i ≤ j), ∀ k : ι, ∀ (hik : i ≤ k),
    ∃ (h : F k →ₗ[R] F j), f i j hij = h.comp (f i k hik)
\end{lean}
Finally, condition (4) connects the module-theoretic definition back to the
inverse system theory of Section~5. For a fixed $R$-module $N$, the family
$(\mathrm{Hom}(M_i, N))_{i \in \iota}$ with transition maps given by precomposition
with $f_{ij}$ forms an inverse system, and condition (4) requires this inverse
system to be Mittag-Leffler in the sense of \texttt{InverseSystem.IsMittagLeffler}
for every $N$.

We formalize the equivalence of all these conditions, following 
\cite[Proposition 6.11]{perry}, and \cite[Tag 0599]{SP}:

\begin{thm}\label{thm:ML-equiv}
Let $R$ be a commutative ring, let $\iota$ be a directed index set, and let
$M = \varinjlim_{i \in \iota} M_i$ be a presentation of an $R$-module $M$ as a
directed colimit of finitely presented modules. Then the following conditions are
equivalent:
\begin{enumerate}
\item $M$ is Mittag-Leffler.
\item The presentation satisfies condition \textup{(2)}.
\item The presentation satisfies condition \textup{(3)}.
\item For every $R$-module $N$, the inverse system $(\mathrm{Hom}(M_i, N))$ 
      is Mittag-Leffler.
\end{enumerate}
\end{thm}

We note that the equivalence of (1) and (2) requires some care beyond what is
stated in \cite[Proposition 6.11]{perry}, and \cite[Tag 0599]{SP}: the argument given there for one
direction is not precise in terms of working with the transition maps of direct limits, and we provide a corrected proof. We formalize
the equivalences as three separate biconditionals:
\begin{lean}
theorem isMittagLeffler_iff_isMittagLeffler' :
    IsMittagLeffler' R F f ↔ IsMittagLeffler R (Module.DirectLimit F f)

theorem isMittagLeffler'_iff_isMittagLeffler'' :
    IsMittagLeffler' R F f ↔ IsMittagLeffler'' R F f

theorem isMittagLeffler''_iff_isMittagLeffler'''
    (hf : ∀ i, f i i (le_refl i) = LinearMap.id)
    (hcomp : ∀ i j k (hij : i ≤ j) (hjk : j ≤ k),
      f i k (hij.trans hjk) = (f j k hjk).comp (f i j hij)) :
    IsMittagLeffler'' R F f ↔ IsMittagLeffler''' R F f hf hcomp
\end{lean}

We also establish the following basic permanence properties. Finitely presented
modules are Mittag-Leffler:
\begin{lean}
theorem Module.finitePresentation_isMittagLeffler (M : Type u) 
    [AddCommGroup M] [Module R M] [Module.FinitePresentation R M] : 
    IsMittagLeffler R M
\end{lean}
and projective modules are Mittag-Leffler:
\begin{lean}
theorem Module.proj_is_Mittag_Leffler {R : Type u} {P : Type u} [CommRing R]
    [AddCommGroup P] [Module R P] [Module.Projective R P] : 
    IsMittagLeffler R P
\end{lean}

Our proof of the latter result is different and easier than the ones given in \cite{SP,perry}. They used another characterization of the Mittag-Leffler property in terms of their behavior with respect to tensoring with direct products.
We are able to bypass this complicated formalization 
by giving a direct simple proof of the facts that free modules are Mittag-Leffler 
and direct summands of Mittag-Leffler modules are Mittag-Leffler,
both of which follow directly from the definition of a Mittag-Leffler module.

The major structural result of this section, which is the key ingredient in the
proof of \cref{thm:projchar}, is the following theorem corresponding to
\cite[Lemma 7.1]{perry} and \cite[Tag 059W]{SP}:

\begin{thm}\label{thm:countable-colimit}
Let $R$ be a commutative ring, let $M = \varinjlim_{i \in \iota} M_i$ be a
presentation of $M$ as a directed colimit of finitely presented modules, and
suppose that $M$ is both Mittag-Leffler and countably generated. Then there
exists a countable directed subset $I' \subseteq \iota$ such that
$M \cong \varinjlim_{i \in I'} M_i$.
\end{thm}

In Lean, this is formalized as:
\begin{lean}
theorem mittagLeffler_countablyGenerated_has_countable_colimit
    {ι : Type u} [Preorder ι] [IsDirected ι (· ≤ ·)] [Nonempty ι] 
    [DecidableEq ι]
    (F : ι → Type u) [∀ i, AddCommGroup (F i)] [∀ i, Module R (F i)]
    [∀ i, Module.FinitePresentation R (F i)]
    (f : ∀ i j, i ≤ j → F i →ₗ[R] F j)
    [DirectedSystem F (fun i j h => f i j h)]
    (hML : IsMittagLeffler R (Module.DirectLimit F f))
    (hCountGen : Module.IsCountablyGenerated R (Module.DirectLimit F f)) :
    ∃ (I' : Set ι) (_ : I'.Countable) (_ : DirectedOn (· ≤ ·) I') 
    (_ : I'.Nonempty),
      Nonempty (Module.DirectLimit F f ≃ₗ[R]
        Module.DirectLimit (fun i : I' => F i.val)
          (fun i j (hij : i ≤ j) => f i.val j.val hij))
\end{lean}

The proof proceeds in two main steps. First, using the countable generation
hypothesis, one finds a countable set of generators of $M$, each of which is
represented at some index in $\iota$, giving an initial countable subset
$I_0 \subseteq \iota$. Second, one enlarges $I_0$ iteratively: at each step,
for each index $i$ in the current set, the Mittag-Leffler condition provides a
stabilization index $j_i \geq i$, which is added to the set. Taking the union
of this sequence yields the required countable directed subset $I'$, and the
identification $M \cong \varinjlim_{i \in I'} M_i$ follows from the fact that
$I'$ was constructed to contain representatives of all generators of $M$.
Contained in it is a separate lemma that we prove,
\texttt{lemma countable\_subset\_directed\_closure},
which says that if you start with a countable subset inside a directed poset,
then you can enlarge it to a countable directed subset
containing the original set.

\section{Projective modules and faithfully flat maps}

The results of this section are in the file basechange.lean in \cite{code}.
We now have all the ingredients to prove the two main theorems. We begin with
the characterization of projective modules, \cref{thm:projchar}, which we
formalize as:

\begin{lean}
theorem proj_iff_Mittag_Leffler {R : Type u} {P : Type u} [CommRing R]
    [AddCommGroup P] [Module R P] : Module.Projective R P ↔
    (Module.Flat R P) ∧ (Module.IsMittagLeffler R P) ∧
    (IsDirectSumOfCountablyGenerated R P)
\end{lean}

The proof of the forward direction is immediate: projectivity implies flatness
by a standard result already in Mathlib, Mittag-Leffler by
\texttt{proj\_is\_Mittag\_Leffler}, and the direct sum decomposition by
\texttt{proj\_is\_dirSumCountable} from Section~1. For the backward direction,
suppose $P = \bigoplus_{i \in I} M_i$ is a direct sum of countably generated
modules, each of which inherits flatness and the Mittag-Leffler property from
$P$ as a direct summand. It therefore suffices to show each $M_i$ is
projective. This reduces the problem to the following key result, which handles
the countably generated case:

\begin{lean}
theorem proj_if_countable_flat_Mittag_Leffler {R : Type u} {P : Type u}
    [CommRing R] [AddCommGroup P] [Module R P]
    (PFlat : Module.Flat R P) (PML : Module.IsMittagLeffler R P)
    (PCount : ∃ s : Set P, s.Countable ∧ span R s = ⊤) :
    Module.Projective R P
\end{lean}

The proof applies Lazard's theorem to write $P$ as a directed colimit of
finitely presented free modules, then uses
\texttt{mittagLeffler\_countablyGenerated\_has\_countable\_colimit} from
Section~7 to replace the index set by a countable one.
The key step is then to apply the projectivity criterion via the lifting
property: given a surjection $f : N_2 \to N_3$ and a map $g : P \to N_3$,
we need to lift $g$ through $f$. Writing $P = \varinjlim_{i \in I'} G_i$ over
a countable directed index set, the maps $g \circ \iota_i : G_i \to N_3$
assemble into an element of the inverse limit $\varprojlim \mathrm{Hom}(G_i,
N_3)$. Since $P$ is Mittag-Leffler, the inverse systems
$(\mathrm{Hom}(G_i, N))$ are Mittag-Leffler for any $N$, by the equivalence
of conditions (1) and (4) in \cref{thm:ML-equiv}. The surjectivity theorem
(Theorem~5.2) then provides a compatible family of lifts in
\[
\varprojlim \mathrm{Hom}(G_i, N_2),
\]
which assembles via
\texttt{DirectLimit.lift} into the required map $P \to N_2$.

We now turn to the main theorem. Two key descent lemmas are needed before
assembling the proof of Theorem~I. The first states that faithful flatness
reflects the Mittag-Leffler property. 
It is an improvement of \cite[Lemma 9.2]{perry},
already present in \cite[tag 05A5]{SP},
and whose proof is attributed there to Juan Pablo Acosta Lopez.

\begin{thm}
Let $R \to S$ be a faithfully flat map of commutative rings, and let $M$ be an
$R$-module. If $S \otimes_R M$ is Mittag-Leffler over $S$, then $M$ is
Mittag-Leffler over $R$.
\end{thm}

which we formalize as:

\begin{lean}
theorem ML_of_faithfullyFlat_tensorProduct {R : Type u} {S : Type u}
    {M : Type u} [CommRing R] [CommRing S] [Algebra R S]
    [Module.FaithfullyFlat R S] [AddCommGroup M] [Module R M]
    (h : Module.IsMittagLeffler S (TensorProduct R S M)) :
    Module.IsMittagLeffler R M
\end{lean}

The proof writes $M$ as a directed colimit $\varinjlim G_i$ of finitely
presented modules and uses the equivalence of conditions (1) and (2) in
\cref{thm:ML-equiv}. The hypothesis gives a stabilization index for each $i$
in the base-changed system $(S \otimes_R G_i)$, and the key step is to
transfer this to the original system using the pushout characterization of
domination \cref{thm:domPushout} together with
\texttt{UniversallyInjective.of\_baseChange\_faithfullyFlat} from Section~4.

The second descent lemma states that faithful flatness reflects countable
generation \cite[tag 0GVD]{SP}:

\begin{thm}
Let $R \to S$ be a faithfully flat map of commutative rings, and let $P$ be an
$R$-module. If $S \otimes_R P$ is countably generated over $S$, then $P$ is
countably generated over $R$.
\end{thm}

which we formalize as:

\begin{lean}
theorem countably_generated_faithfully_flat {R : Type u} {S : Type u}
    {P : Type u} [CommRing R] [CommRing S] [Algebra R S]
    [Module.FaithfullyFlat R S] [AddCommGroup P] [Module R P]
    (PSIsCountable : ∃ s : Set (TensorProduct R S P),
      s.Countable ∧ span S s = ⊤) :
    ∃ t : Set P, t.Countable ∧ span R t = ⊤
\end{lean}

The elementary proof enumerates a countable generating set of $S \otimes_R P$, expresses
each generator as a finite sum of pure tensors, collects all the $P$-components
into a countable subset $t \subseteq P$, and then uses faithful flatness to
show that $t$ generates $P$ over $R$.

Combining these two descent lemmas with
\texttt{proj\_if\_countable\_flat\_Mittag\_Leffler} immediately gives the
countably generated case of Theorem~I:

\begin{lean}
theorem proj_countable_faithfully_flat {R : Type u} {S : Type u}
    {P : Type u} [CommRing R] [CommRing S] [Algebra R S]
    [Module.FaithfullyFlat R S] [AddCommGroup P] [Module R P]
    (PIsCountable : ∃ s : Set (TensorProduct R S P),
      s.Countable ∧ span S s = ⊤)
    (baseChangeProj : Module.Projective S (TensorProduct R S P)) :
    Module.Projective R P
\end{lean}

The general case of Theorem~I requires a transfinite induction argument to
reduce to the countably generated case. The key idea is to build a
Kaplansky-style filtration $0 = P_0 \subseteq P_1 \subseteq \cdots \subseteq
P_\gamma = P$ indexed by ordinals up to $\gamma = \mathrm{succ}(\mathrm{type}(P))$,
with the property that each successive quotient $P_{\alpha+1}/P_\alpha$ is
countably generated and projective over $R$. This filtration is constructed by
transfinite recursion: at each successor step, a well-ordering of $P$ is used
to pick the least element not yet in $P_\alpha$, and
\texttt{countably\_generated\_lift} extends the current stage by a countably
generated piece while maintaining the correspondence between $R$-submodules of
$P$ and sub-direct-sums of the decomposition $S \otimes_R P = \bigoplus_{i \in
I} Q_i$. At limit stages, one takes the union. Each successive quotient is
shown to be projective using \texttt{proj\_countable\_faithfully\_flat},
relying on the isomorphism $S \otimes_R (P_{\alpha+1}/P_\alpha) \cong
({\bigoplus}_{i \in I'} Q_i) / ({\bigoplus}_{i \in I} Q_i)$ for appropriate
index sets $I \subseteq I'$, which is itself a direct summand of $S \otimes_R
P$ and hence projective. The filtration therefore constitutes a Kaplansky
devissage for $P$, which by \texttt{kaplansky\_devissage\_iff\_direct\_sum}
expresses $P$ as a direct sum of countably generated projective modules, and
thus projective. We thus obtain
\cref{thm:main}:

\begin{lean}
theorem proj_faithfully_flat {R : Type u} {S : Type u} {P : Type u}
    [CommRing R] [CommRing S] [Algebra R S] [Module.FaithfullyFlat R S]
    [AddCommGroup P] [Module R P] :
    Module.Projective R P ↔ Module.Projective S (TensorProduct R S P)
\end{lean}

\bibliographystyle{amsplain}
\bibliography{main}

\end{document}